\documentclass[12pt,reqno,oneside,titlepage]{amsart}

\newtheorem{theorem}{Theorem}[section]
\newtheorem{proposition}[theorem]{Proposition}
\newtheorem{definition}[theorem]{Definition}

\newtheorem{corollary}[theorem]{Corollary}
\newtheorem{conjecture}[theorem]{Conjecture}

\newcommand{\bdf}{\mathbf{f}}
\newcommand{\bdg}{\mathbf{g}}

\newcommand\lcb{\left\{}
\newcommand\rcb{\right\}}

\newcommand\lp{\left(}
\newcommand\rp{\right)}

\newcommand\cmp{\mathbf{p}}
\newcommand\cmpq{\mathbf{q}}

\newcommand\Res{\mathop{\mathrm{Res}}}

\newcommand\sth{^{\text{th}}}

\newcommand\natnums{\mathbb{N}}

\newcommand\cD{\mathcal{D}}

\title[Composition sum identities]{Composition sum identities related to the
  distribution of coordinate values in a discrete simplex.}

\author{R. Milson\\Dept. Mathematics \&  Statistics\\Dalhousie University}


\address{Dept of Math., Dalhousie U., Halifax, Canada, B3H 3J5}

\email{milson@mathstat.dal.ca}

\thanks{Email: {\tt milson@mathstat.dal.ca}}
\thanks{This research supported by a Dalhousie University startup grant.}

\begin{document}

\begin{abstract}
  Utilizing spectral residues of parameterized, recursively defined
  sequences, we develop a general method for generating identities of
  composition sums. Specific results are obtained by focusing on
  coefficient sequences of solutions of first and second order,
  ordinary, linear differential equations.
  
  Regarding the first class, the corresponding identities amount to a
  proof of the exponential formula of labelled counting.  The
  identities in the second class can be used to establish certain
  geometric properties of the simplex of bounded, ordered, integer
  tuples.
  
  We present three theorems that support the conclusion that the inner
  dimensions of such an order simplex are, in a certain sense, more
  ample than the outer dimensions.  As well, we give an algebraic
  proof of a bijection between two families of subsets in the order
  simplex, and inquire as to the possibility of establishing this
  bijection by combinatorial, rather than by algebraic methods.
\end{abstract}
\maketitle

\pagestyle{myheadings}

\section{Introduction}

The present paper is a discussion of composition sum identities that
may be obtained by utilizing spectral residues of parameterized,
recursively defined sequences.  Here we are using the term
``composition sum'' to refer to a sum whose index runs over all
ordered lists of positive integers $p_1, p_2, \ldots, p_l$ that such
that for a fixed $n$,
$$p_1+\ldots + p_l=n.$$
Spectral residues will be discussed in detail
below.

Compositions sums are a useful device, and composition sum identities
are frequently encountered in combinatorics.  For example the Stirling
numbers (of both kinds) have a natural representation by means of such sums:
\cite[\S 51, \S 60]{Jordan}:
$$
s^l_n = \frac{n!}{l!}
\sum_{p_1+\ldots+p_l=n}
\frac{1}{p_1\,p_2\,\ldots p_l}; \qquad
\mathfrak{S}^l_n = \frac{n!}{l!}
\sum_{p_1+\ldots+p_l=n} 
\frac{1}{p_1!\,p_2!\,\ldots p_l!}.  
$$
There are numerous other examples.  In general, it is natural to
use a composition sum to represent the value of quantities $f_n$ that
depend in a linearly recursive manner on quantities $f_1, f_2, \ldots,
f_{n-1}$.  By way of illustration, let us mention that this point of
view leads immediately to the interpretation of the $n\sth$ Fibonacci
number as the cardinality of the set of compositions of $n$ by
$\{1,2\}$ \cite[2.2.23]{Goulden}

To date, there are few systematic investigations of composition sum
identities.  The references known to the present author are
\cite{Hoggatt} \cite{Homenko} \cite{Moser}; all of these papers obtain
their results through the use of generating functions.  In this
article we propose a new technique based on spectral residues, and
apply this method to derive some results of an enumerative nature.
Let us begin by describing one of these results, and then pass to a
discussion of spectral residues.

Let $S^3(n)$ denote the discrete simplex of bounded, ordered triples
of natural numbers:
$$S^3(n) = \{(x,y,z)\in \natnums^3: 0\leq x<y<z\leq n\}.$$
In regard
to this simplex, we may inquire as to what is more probable: a
selection of points with distinct $y$ coordinates, or a selection of
points with distinct $x$ coordinates.  The answer is given by the
following.
\begin{theorem}
\label{thrm:id1}
For every cardinality $l$ between $2$ and $n-1$, there are more
$l$-element subsets of $S^3(n)$ with distinct $y$ coordinates, than
there are $l$-element subsets with distinct $x$ coordinates.
\end{theorem}

Let us consider this result from the point of view of generating
functions.  The number of points with $y=j$ is $j(n-j)$.  Hence the
generating function for subsets with distinct $y$-values is
$$Y(t)=\prod_{j=1}^{n-1} (1+j(n-j)t),$$
where $t$ counts the selected
points.  The number of points with $x=n-j$ is $j(j-1)/2$.  Hence, the
generating function for subsets with distinct $x$-values is
$$X(t) = \prod_{j=2}^n \lp 1+\frac{j(j-1)}2\, t\rp.$$
The above
theorem is equivalent to the assertion that the coefficients of $Y(t)$
are greater than the coefficients of $X(t)$.  The challenge is to find
a way to compare these coefficients.

We will see below this can be accomplished by re-expressing the
coefficients in question as composition sums, and then employing a
certain composition sum identity to make the comparison.  We therefore
begin by introducing a method for systematically generating such
identities.

\section{The method of spectral residues}
Let us consider a sequence of quantities $f_n$, recursively defined by
\begin{gather}
\label{eq:recrel}
f_0 = 1,\qquad
(\nu-n)f_n = \sum_{j=0}^{n-1} a_{jn} f_j,\quad n=1,2,\ldots 
\end{gather}
where the $a_{jk},\, 0\leq j<k$ is a given array of constants, and
$\nu$ is a parameter.  The presence of the parameter has some
interesting consequences.

For instance, it is evident that if $\nu$ is a natural number, then
there is a possibility that the relations \eqref{eq:recrel} will not
admit a solution.  To deal with this complication we introduce the
quantities
$$\rho_n = \Res(f_n(\nu),\nu=n),$$
and henceforth refer to them as
spectral residues.  The list $\rho_1,\rho_2,\ldots$ will be called the
spectral residue sequence.
\begin{proposition}
  If $\nu=n$ then the relations \eqref{eq:recrel} do not admit a solution if
  $\rho_n\neq 0$, and admit multiple solutions if $\rho_n=0$.
\end{proposition}
\begin{proof}
  If $\nu=n$, the relations in question admit a solution if and only
  if
  $$\sum_{j=1}^{n-1} a_{jn} f_j \Big|_{\nu=n} = 0.$$
  The left-hand
  side of the above equation is precisely, $\rho_n$, the $n\sth$
  spectral residue.  It follows that if $\rho_n=0$, then the value of
  $f_n$ can be freely chosen, and that the solutions are uniquely
  determined by this value.
\end{proof}

The above proposition is meant to indicate how spectral residues arise
naturally in the context of parameterized, recursively defined
sequences.  However, our interest in spectral residues is motivated by
the fact that they can be expressed as composition sums.  To that end,
let $\cmp=(p_1,\ldots,p_l)$ be an ordered list of natural numbers. We
let
$$s_j=p_1+\ldots+p_j,\quad j=1,\ldots, l$$
denote the $j\sth$ left partial sum and set
$$|\cmp|=s_l=p_1+\ldots+p_l.$$
Let us also define the following abbreviations:
$$
s_\cmp = \prod_{j=1}^{l-1} s_j, \qquad
a_\cmp = \prod_{j=1}^{l-1} a_{s_js_{j+1}}
$$
\begin{proposition}
$$\rho_n = \sum_{|\cmp|=n} a_\cmp / s_\cmp.$$  
\end{proposition}
\noindent
Composition sum identities arise in this setting because spectral
residue sequences enjoy a certain invariance property.

Let $\bdf=(f_1,f_2,\ldots)$  and $\bdg=(g_1,g_2,\ldots)$ be sequences defined,
respectively by relation \eqref{eq:recrel} and by
\begin{gather*}
g_0 = 1,\qquad
(\nu-n)g_n = \sum_{j=0}^{n-1} b_{jn} g_j,\quad n=1,2,\ldots \\
\end{gather*}
\begin{definition}
\label{def:unirel}
We will say that $\bdf$ and $\bdg$ are unipotently equivalent if $g_n
= f_n$ plus a $\nu$-independent linear combination of $f_1,\ldots,
f_{n-1}$.
\end{definition}
The motivation for this terminology is as follows.  It is natural to
represent the coefficients $a_{ij}$ and $b_{ij}$ by infinite, lower
nilpotent matrices, call them $A$ and $B$.  Let $D_\nu$ denote the
diagonal matrix with entry $\nu-n$ in position $n+1$. The sequences
$\bdf$ and $\bdg$ are then nothing but generators of the kernels of
$D_\nu-A$ and $D_\nu-B$, respectively, The condition that $\bdf$ and
$\bdg$ are unipotently equivalent amounts to the condition that
$D_\nu-A$ and $D_\nu-B$ are related by a unipotent matrix factor.

Unipotent equivalence is, evidently, an equivalence relation on the
set of sequences of type \eqref{eq:recrel}.  
\begin{proposition}
  The spectral residue sequence is an invariant of the corresponding
  equivalence classes.
\end{proposition}
\begin{proof}
  The recursive nature of the $f_k$ ensures that $\Res(f_k;\nu=n)$
  vanishes for all $k<n$.  The proposition now follows by inspection
  of Definition \ref{def:unirel}.
\end{proof}
\noindent The application of this  result to composition identities is
  immediate. 
\begin{corollary}
  If $a_{ij}$ and $b_{ij}$ are nilpotent arrays of constants such that the
  corresponding $\bdf$ and $\bdg$ are unipotently equivalent, then
  necessarily
  $$\sum_{|\cmp|=n} a_\cmp / s_\cmp = \sum_{|\cmp|=n} b_\cmp / s_\cmp.$$
\end{corollary}

Due to its general nature, the above result does not, by itself, lead
to interesting composition sum identities.  In the search for useful
applications we will limit our attention to recursively defined
sequences arising from series solutions of linear differential
equations.  Consideration of both first and second order equations in
one independent variable will prove fruitful.  Indeed, in the next
section we will show that the first-order case naturally leads to the
exponential formula of labelled counting \cite[\S 3]{wilf}.  The
second-order case will be considered after that; it leads naturally to
the type of result discussed in the introduction.

\section{Spectral residues of first-order equations.}
Let $U=U_1z+U_2z^2+\ldots$ be a formal power series with zero constant
term, and let $\phi(z)$ be the series solution of the following
parameterized, first-order, differential equation:
$$z\phi'(z)+[U(z)-\nu]\phi(z) +\nu = 0,\quad \phi(0)=1.$$
Equivalently, the coefficients of $\phi(z)$ must satisfy
$$\phi_0=1,\qquad (\nu-n)\phi_n = \sum_{j=0}^{n-1}
U_{n-j}\phi_j.$$

In order to obtain a composition sum identity we seek a related
equation whose solution will be unipotently related to $\phi(z)$.  It
is well known that a linear, first-order differential equation can be
integrated by means of a gauge transformation.  Indeed, setting
\begin{align*}
  \sigma(z) &= \sum_{k=1}^\infty U_k\frac{z^k}{k},\\
  \psi(z) &= \exp(\sigma(z)) \phi(z)
\end{align*}
our differential equation is transformed into
$$z\psi'(z) - \nu\psi(z) + \nu \exp(\sigma(z))=0.$$
Evidently, the
coefficients of $\phi$ and $\psi$ are unipotently related, and hence
we obtain the following composition sum identity.
\begin{proposition}
  Setting $U_\cmp=\prod_i U_{p_i}$ for $\cmp=(p_1,\ldots,p_l)$ we have
  \begin{equation}
    \label{eq:1ordid}
    \sum_n \sum_{|\cmp|=n} \frac{U_\cmp}{s_\cmp } \frac{z^n}n=
  \exp\lp\sum_k U_k\frac{z^k}k\rp .    
  \end{equation}
\end{proposition}

The above identity has an interesting interpretation in the context of
labelled counting, e.g. the enumeration of labelled graphs.  In our
discussion we will adopt the terminology introduced in H. Wilf's book
\cite{wilf}.  For each natural number $k\geq 1$ let $\cD_k$ be a set
--- we will call it a deck --- whose elements we will refer to as {\em
  pictures of weight $k$}.  A {\em card of weight $k$} is a pair
consisting of a picture of weight $k$ and a $k$-element subset of
$\natnums$ that we will call {\em the label set} of the card.  A {\em
  hand of weight $n$ and size $l$} is a set of $l$ cards whose weights
add up to $n$ and whose label sets form a partition of
$\{1,2,\ldots,n\}$ into $l$ disjoint groups.  The goal of labelled
counting is to establish a relation between the cardinality of the
sets of hands and the cardinality of the decks.

For example, when dealing with labelled graphs, $\cD_k$ is the set of
all connected $k$-graphs whose vertices are labelled by $1, 2, \ldots,
k$.  A card of weight $k$ is a connected $k$-graph labelled by any $k$
natural numbers.  Equivalently, a card can be specified as a picture
and a set of natural number labels. To construct the card we label
vertex $1$ in the picture by the smallest label, vertex $2$ by the
next smallest label, etc.  Finally, a hand of weight $n$ is an
$n$-graph (not necessarily connected) whose vertices are labelled by
$1,2,\ldots,n$.

Let $d_k$ denote the cardinality of $\cD_k$ and  set
$$d(z)=\sum_k d_k\frac{z^k}{k!}$$
Similarly let
$h_{nl}$ denote the cardinality of the set of hands of weight $n$ and
size $l$, and set
$$h(y,z)=\sum_{nl} h_{nl}\,y^l\,\frac{z^n}{n!}.$$
The exponential
formula of labelled counting is an identity that relates the above
generating functions.  Here it is:
\begin{equation}
  \label{eq:expform}
  h(y,z) = \exp(y\,d(z)).  
\end{equation}

To establish the equivalence of \eqref{eq:1ordid} and
\eqref{eq:expform} we need to introduce some extra terminology.
Consider a list of $l$ cards with weights $p_1,\ldots,p_l$ and label
sets $S_1, \ldots, S_l$.  We will say that such a list forms an {\em
ordered hand} if
$$\min(S_i)<\min(S_{i+1}),\quad \mbox{ for all } i=1,\ldots,l-1.$$
Evidently, each hand (a set of cards) corresponds to a unique ordered
hand (an ordered list of the same cards), and hence we seek a way to
enumerate the set of all ordered hands of weight $n$ and size $l$.

Let us fix a composition $\cmp=(p_1,\ldots,p_l)$ of a natural number
$n$, and consider a permutation $\pi=(\pi_1,\ldots,\pi_n)$ of
$\{1,\ldots,n\}$.  Let us sort $\pi$ according to the following
scheme.  Exchange $\pi_1$ and $1$ and then sort
$\pi_2,\ldots,\pi_{p_1}$ into ascending order.  Next exchange
$\pi_{p_1+1}$ and the minimum of $\pi_{p_1+1},\ldots,\pi_n$ and then
sort $\pi_{p_1+2},\ldots, \pi_{p_2}$ into ascending order.  Continue
in an analogous fashion $l-2$ more times.  The resulting permutation
will describe a division of $\{1,\ldots,n\}$ into $l$ ordered blocks,
with the blocks themselves being ordered according to their smallest
elements.  Call such a permutation $\cmp$-ordered.  Evidently, each
$\cmp$-ordered permutation can be obtained by sorting
$$s_\cmp\times n\times \prod_i (p_i-1)!$$
different permutations. 

Next, let us note that an ordered hand can be specified in terms of
the following ingredients: a composition $\cmp$ of $n$,  one of
$\prod_i d_{p_i}$ choices of pictures of weights $p_1,\ldots,p_l$, and
a $\cmp$-ordered permutation. It follows that
$$h_{nl} = \sum_{\stackrel{\scriptstyle|\cmp|=n}{
    \cmp=(p_1,\ldots,p_l)}} \frac{n!}{s_\cmp\times n\times\prod_i
  (p_i-1)!} \, \prod_i d_{p_i}.$$
Finally, we can establish the
equivalence of \eqref{eq:1ordid} and \eqref{eq:expform} by setting
$$U_k = \frac{d_k}{(k-1)!}\,y.$$

\section{Spectral residues of second-order equations.}
Let $U=U_1z+U_2z^2+\ldots$ be a formal power series with zero constant
term, and let $\phi(z)$ be the series solution of the following
second-order, linear differential equation:
\begin{equation}
  \label{eq:soeqn}
  z^2\phi''(z)+(1-\nu)z\phi'z + U(z)\phi(z)= 0,\quad \phi(0)=1.
\end{equation}
Equivalently, the coefficients of $\phi(z)$ are determined by
$$\phi_0=1,\qquad n(\nu-n)\phi_n = \sum_{j=0}^{n-1} U_{n-j}\phi_j.$$

Two remarks are in order at this point. First, the class of equations
described by \eqref{eq:soeqn} is closely related to the class of
self-adjoint second-order equations.  Indeed, conjugation by a gauge
factor $z^{\nu/2}$ transforms \eqref{eq:soeqn} into self-adjoint form
with potential $U(z)$ and energy $\nu^2/4$.  The solutions of the
self-adjoint form are formal series multiplied by $z^{\nu/2}$, so
nothing is lost by working with the ``nearly'' self-adjoint form
\eqref{eq:soeqn}.

Second, there is no loss of generality in restricting our focus to the
self-adjoint equations.  Every second-order linear equation can be
gauge-transformed into self-adjoint form, and as we saw above,
spectral residue sequences are invariant with respect to gauge
transformations.  Indeed, as we shall demonstrate shortly, the
potential $U(z)$ is uniquely determined by its corresponding residue
sequence.
\begin{proposition}
\label{prop:2specres}
  The spectral residues corresponding to \eqref{eq:soeqn} are
  $$\rho_n = \frac1n \sum_{|\cmp|=n}
  \frac{U_\cmp}{s_\cmp\,s_{\cmp'}},$$
  where as before, for
  $\cmp=(p_1,\ldots,p_l)$, we write $U_\cmp$ for $\prod_i U_{p_i}$,
  and write $\cmp'$ for the reversed composition
  $(p_l,p_{l-1},\ldots,p_1)$.
\end{proposition}
Since $\rho_n = U_n/n$ plus a polynomial of $U_1,\ldots,U_{n-1}$, it
is evident that the spectral residue sequence completely determines
the potential $U(z)$.  An explicit formula for the inverse relation is
given in \cite{milson2}.

Interesting composition sum identities will appear in the present
context when we consider exactly-solvable differential equations.  We
present three such examples below, and discuss the enumerative
interpretations in the next section.  In each case the exact
solvability comes about because the equation is gauge-equivalent to
either the hypergeometric, or the confluent hypergeometric equation.
Let us also remark --- see \cite{milson2} for the details --- that
these equations occupy an important place within the canon of
classical quantum mechanics, where they correspond to various
well-known exactly solvable one-dimensional models.

\begin{proposition}
  \label{prop:id1}
  $$\sum_{\stackrel{\scriptstyle\cmp=(p_1,\ldots,p_l)}{\strut
      |\cmp|=n}} \!\!\!\!\
      \frac{(n-1)!}{s_\cmp}\,\frac{(n-1)!}{s_{\cmp'}}\lp\prod_i
      p_i\rp t^l = \prod_{j=1}^{n} \lcb t+j(j-1) \rcb$$ 
\end{proposition}
\begin{proof}
  By Proposition \ref{prop:2specres}, the left hand side of the above
  identity is $n!(n-1)!$ times the $n\sth$ spectral residue
  corresponding to the potential
  $$U(z) = \frac{tz}{(z-1)^2}=t\sum_k kz^k.$$
  Setting
  $$t=\alpha(1-\alpha)$$
  and making a change of gauge
  $$\phi(z)=(1-z)^\alpha \psi(z)$$
  transforms \eqref{eq:soeqn} into
  $$z^2\psi''(z) + (1-\nu)\psi'(z) - \frac{z}{1-z}\lcb2\alpha z\,\psi'(z) +
  \alpha(\alpha-\nu)\,\psi(z)\rcb=0.$$
  Multiplying through by $(1-z)/z$ and setting 
  $$\gamma=1-\nu,\quad \beta=\alpha-\nu,$$
  we recover the usual hypergeometric equation
  $$z(1-z)\,\psi''(z) + \lcb \gamma+(1-\alpha-\beta)z \rcb \psi'(z)
  -\alpha\beta\,\psi(z)=0. $$
  It follows that
  $$\psi_n = \frac{(\alpha)_n (\alpha-\nu)_n}{n! (1-\nu)_n},$$
  and hence the $n\sth$ spectral residue is given by
  $$\rho_n=(-1)^n\frac{\prod_{j=1}^n  (\alpha-j)(\alpha+j-1)}{n! (n-1)!},$$
  or equivalently by
  $$\rho_n = \frac{\prod_{j=1}^n (t+j(j-1))}{n! (n-1)!}.$$
  The
  asserted identity now follows from the fundamental invariance
  property of spectral residues.
\end{proof}

\begin{proposition}
  \label{prop:id2}
  $$\sum_{\stackrel{\scriptstyle
      \cmp=(p_1,\ldots,p_l)}{\stackrel{\scriptstyle \strut
        p_i\in\{1,2\}}{|\cmp|=n}}} \!\!\!\!\!
  \frac{(n-1)!}{s_\cmp}\,\frac{(n-1)!}{s_{\cmp'}}\, t^{n-l} = \prod_k
  ( 1+k^2t),$$
  where the right hand index $k$ varies over all positive integers
      $n-1, n-3, n-5,\ldots$. 
\end{proposition}
\begin{proof}
  As in the preceding proof, Proposition \ref{prop:2specres} shows
  that the left hand side of the present identity is $n!(n-1)!$ times
  the $n\sth$ spectral residue corresponding to the potential
  $$U(z) = z+tz^2.$$
  Setting
  $$t=-\omega^2,$$
  and making a change of gauge
  $$\phi(z) = \exp(\omega z) \psi(z)$$
  transforms \eqref{eq:soeqn} into
  $$z^2\psi''(z) + (1-\nu)z\psi'(z) + 2\omega z^2\psi'(z) + z\,(\omega
  (1-\nu) +1)\,\psi(z)=0.$$
  Dividing through by $z$ and setting
  $$\gamma=1-\nu,\quad 1=\omega(2\alpha+\nu-1),$$
  we obtain the following scaled variation of the confluent
  hypergeometric equation:
  $$z\psi''(z) + (\gamma+2\omega z)\psi'(z) +
  2\omega\alpha\,\psi(z)=0.$$
  It follows that 
  $$\psi_n = \frac{(-2\omega)^2 (\alpha)_n}{n!(\gamma)_n},$$
  and hence that
  \begin{align*}
    \rho_n &= \frac{\prod_{k=0}^{n-1} (1+\omega(2k+1-n))}{n!(n-1)!}\\
    &= \frac{\prod_{k=0}^{\lfloor\frac{n-1}2\rfloor}
      (1+t(n-1-2k)^2)}{n!(n-1)!}.
  \end{align*}
  The
  asserted identity now follows from the fundamental invariance
  property of spectral residues.
\end{proof}

\begin{proposition}
  \label{prop:id3}
  $$\sum_{\stackrel{\scriptstyle
      \cmp=(p_1,\ldots,p_l)}{\stackrel{\scriptstyle \strut
        p_i\; \mbox{\tiny\rm odd}}{|\cmp|=n}}} \!\!\!\!\!
  \frac{(n-1)!}{s_\cmp}\,\frac{(n-1)!}{s_{\cmp'}}\lp\prod_i
      p_i\rp t^{\frac{n-l}2} = \prod_k
  \lcb 1+(k^4-k^2)t\rcb,$$
  where the right hand index $k$ ranges over all positive integers
      $n-1, n-3, n-5, \ldots$.
\end{proposition}
\begin{proof}
  By Proposition \ref{prop:2specres}, the left hand side of the
  present identity is $n!(n-1)!\,t^{n/2}$ times the $n\sth$ spectral
  residue corresponding to the potential
  $$U(z) = \frac1{2\sqrt{t}}\lp\frac{z}{(1-z)^2}+\frac{z}{(1+z)^2}\rp=
  \frac{1}{\sqrt{t}}\sum_{k\;\mbox{\tiny \rm odd}} kz^k.$$
  The rest of
  the proof is similar to, but somewhat more involved than the proofs
  of the preceding two Propositions.  Suffice it to say that with the
  above potential, equation \eqref{eq:soeqn} can be integrated by
  means of a hypergeometric function.  This fact, in turn, serves to
  establish the identity in question.  The details of this argument
  are to be found in \cite{milson2}.
\end{proof}
\section{Distribution of coordinate values in a discrete simplex}
In this section we consider enumerative interpretations of the
composition sum identities derived in Proposition \ref{prop:id1},
\ref{prop:id2}, \ref{prop:id3}.  Let us begin with some general
remarks about compositions.

There is a natural bijective correspondence between the set of
compositions of $n$ and the powerset of $\{1,\ldots,n-1\}$.  The
correspondence works by mapping a composition $\cmp=(p_1,\ldots,p_l)$
to the set of left partial sums $\{s_1,\ldots,s_{l-1}\}$, henceforth
to be denoted by $L_\cmp$.  It may be useful to visualize this
correspondence it terms of a ``walk'' from $0$ to $n$: the composition
specifies a sequence of displacements, and $L_\cmp$ is the set of
points visited along the way.  One final item of terminology: we will
call two compositions $\cmp$, $\cmpq$ of $n$ complimentary, whenever
$L_\cmp$ and $L_\cmpq$ disjointly partition $\{1,\ldots,n-1\}$.

Now let us turn to the proof of Theorem \ref{thrm:id1}.
As was mentioned in the introduction, this Theorem is equivalent to
the assertion that the coefficients of 
$$
Y(t)=\prod_{j=1}^{n-1} (1+j(n-j)t)
$$
are greater than the corresponding coefficients of
$$X(t) = \prod_{j=2}^n \lp 1+\frac{j(j-1)}2\, t\rp.$$
Rewriting the former function as a composition sum we have
$$Y(t)=\!\!\!\!\!\sum_{\stackrel{\scriptstyle\cmp=(p_1,\ldots,p_l)}{\strut
    |\cmp|=n}} \!\!\!\!\! s_\cmp\,s_{\cmp'}\,t^l,$$
or equivalently 
$$
Y(t)=\!\!\!\!\!\sum_{\stackrel{\scriptstyle\cmp=(p_1,\ldots,p_l)}{\strut
    |\cmp|=n}} \!\!\!\!\!
\frac{(n-1)!}{s_\cmp}\frac{(n-1)!}{s_{\cmp'}}\,t^{n-l}.$$
On the other
hand, Proposition \ref{prop:id1} allows us to write
$$X(t)=
  \!\!\!\!\!\sum_{\stackrel{\scriptstyle\cmp=(p_1,\ldots,p_l)}{\strut
      |\cmp|=n}} \!\!\!\!\!
      \frac{(n-1)!}{s_\cmp}\,\frac{(n-1)!}{s_{\cmp'}}\lp\prod_i
      \frac{p_i}{2^{p_i-1}} \rp t^{n-l}.$$
It now becomes a straightforward matter to  compare the coefficients
  of $Y(t)$ to those of $X(t)$.  Indeed the desired conclusion follows
  from the rather obvious inequality:
  $$k\leq 2^{k-1},\quad k=1,2,3\ldots ,$$
  the inequality being strict for $k\geq 3$.

Let us now turn to an enumerative interpretation of the composition
sum identity featured in Proposition \ref{prop:id2}.  In order to
state the upcoming result we need to define two notions of sparseness
for subsets of $S^3(n)$.  Let us call a multiset $M$ of integers {\em
  sparse} if $M$ does not contain duplicates, and if
$$|a-b|\geq 2$$
for all distinct $a,b\in M$.  Let us also say that a
multiset $M$ is {\em $2$-sparse} if $M$ does not contain duplicates,
and if there do not exist distinct $a,b\in M$ such that
$$\lfloor a/2 \rfloor = \lfloor b/2 \rfloor.$$  It isn't hard to see
that sparseness is a more restrictive notion than $2$-sparseness,
i.e. if $M$ is sparse, then it is necessarily $2$-sparse, but not the
other way around.  For example, the set
$$\{1,3,4,7\}$$
is not sparse, but it is $2$-sparse.

We require one other item of notation.  For $A\subset S^3(n)$ we let
$\pi_x(A)$ denote the multiset of $x$-coordinates of points in $A$,
and let $\pi_y(A)$ denote the multiset of $y$-coordinates.  We are now
ready to state
\begin{theorem}
\label{thrm:id2}
For every cardinality $l$ between $2$ and $n-1$, there are more
$l$-element subsets $A$ of $S^3(n)$ such that $\pi_y(A)$ is sparse,
than there are $l$-element subsets $A$ such that $\pi_x(A)$ is sparse.
Indeed, the number of $l$-element subsets $A$ of $S^3(n)$ such that
$\pi_y(A)$ is sparse is equal to the number of $l$-element subsets $A$
of $S^3(n)$ such that $\pi_x(A)$ is merely $2$-sparse.
\end{theorem}
\begin{proof}
  Let $\cmp$ be a composition of $n$. Let us begin by noting that the
  corresponding $L_\cmp$ is sparse if and only if the complimentary
  composition consists of $1$'s and $2$'s only.  It therefore follows
  that the enumerating function for $A\subset S^3(n)$ such that
  $\pi_y(A)$ is sparse is
  $$\sum_{\stackrel{\scriptstyle
      \cmp=(p_1,\ldots,p_l)}{\stackrel{\scriptstyle \strut
        p_i\in\{1,2\}}{|\cmp|=n}}} \!\!\!\!\!
  \frac{(n-1)!}{s_\cmp}\,\frac{(n-1)!}{s_{\cmp'}}\, t^{n-l}.$$
  On the
  other hand the number of $(x,y,z)\in S^3(n)$ such that $x\in
  \{2k,2k+1\}$ for any given $k$ is precisely
  $$\binom{n-2k}{2}+\binom{n-2k-1}{2}= (n-2k-1)^2.$$
  Hence the enumerating function for $A\subset S^3(n)$ such that $\pi_x(A)$
  is $2$-sparse is
  $$
  \prod_{k=0}^{\lfloor (n-1)/2\rfloor} \lp 1+(n-2k-1)^2t\rp.$$
  The
  two enumerating functions are equal by Proposition \ref{prop:id2}.
\end{proof}

Finally, let us consider an enumerative interpretation of the
composition sum identity featured in Proposition \ref{prop:id2}.  The
setting for this result will be $S^5(n)$, the discrete simplex of all
bounded, ordered $5$-tuples $(x_1,x_2,x_3,x_4,x_5)$.  For $A\subset
S^5(n)$ we will use $\pi_i(A),\, i=1,\ldots,5$ to denote the
corresponding multiset of $x_i$ coordinate values.
\begin{theorem}
\label{thrm:id3}
For every cardinality $l$ between $2$ and $n-3$, there are more
$l$-element subsets $A$ of $S^5(n)$ such that $\pi_3(A)$ is sparse,
than there are $l$-element subsets $A$ such that $\pi_1(A)$ is $2$-sparse.
\end{theorem}
\begin{proof}
  Let us note that
  the number of points in $S^5(n)$ such that $x_3=j+1$ is given by
  $$\binom{j+1}{2}\binom{n-j-1}{2}.$$
  Hence, the enumerating function for
  the first class of subsets is given by
  $$X_3(t) = \!\!\!\!\!\sum_{\stackrel{\scriptstyle
      \cmp=(p_1,\ldots,p_l)}{\stackrel{\scriptstyle \strut
        p_i\in\{1,2\}}{|\cmp|=n-2}}} \left\{ \prod_{j\not\in
      L_\cmp}\frac{j(j+1)(n-j-1)(n-j-2)}4\right\} t^{n-2-l}.
   $$
   Now there is a natural bijection between the set of compositions
   of $n-2$ by $\{1,2\}$ and the set of compositions of $n-1$ by odd
   numbers.  The bijection works by prepending a $1$ to a composition
   of the former type, and then performing substitutions of the form
   $$(\ldots, k,2,\ldots) \mapsto (\ldots, k+2,\ldots),\quad k\mbox
   { \rm odd}.$$
   Consequently, we can write
   \begin{equation}
     \label{eq:x3sum}
     X_3(t) = \!\!\!\!\!\sum_{\stackrel{\scriptstyle
      \cmp=(p_1,\ldots,p_l)}{\stackrel{\scriptstyle \strut
        p_i\mbox{ \tiny\rm odd}}{|\cmp|=n-1}}}
   \frac{n!}{s_\cmp}\;\frac{n!}{s_{\cmp'}}\; \lp\frac{t}{4}\rp^{(n-1-l)/2}.
   \end{equation}
   Turning to the other class of subsets, the number of points
   $(x_1,\ldots, x_5)$ that satisfy 
   $$x_1\in\{2j,2j+1\}$$
   is given by
   $$\binom{n-2j}{4}+\binom{n-2j-1}{4} = \frac{(n-2j-2)^4 -
   (n-2j-2)^2}{12}.$$ 
   Consequently the enumerating function for subsets $A$ such that
   $\pi_1(A)$ is $2$-sparse is given by
   $$X_1(t) = \prod_k \lp 1+(k^4-k^2)\frac{t}{12}\rp,$$
   where $k$ ranges over all positive integers $n-2, n-4, \ldots.$
   Next, using the identity in Proposition \ref{prop:id3} we have
   $$X_1(t) = \!\!\!\!\!\sum_{\stackrel{\scriptstyle
       \cmp=(p_1,\ldots,p_l)}{\stackrel{\scriptstyle \strut p_i\;
         \mbox{\tiny\rm odd}}{|\cmp|=n-1}}} \!\!\!\!\!
   \frac{(n-1)!}{s_\cmp}\,\frac{(n-1)!}{s_{\cmp'}} \lp \prod_i
   \frac{p_i}{3^{(p_i-1)/2}} \rp \lp\frac{t}{4}\rp^{(n-1-l)/2}.$$
   Using \eqref{eq:x3sum} it now becomes a straightforward matter to
   compare $X_1(t)$ to $X_3(t)$.  Indeed, the desired conclusion
   follows from the following evident inequality:
   $$k\leq 3^{(k-1)/2},\quad k=1,3,5,\ldots ,$$
   the inequality being strict for $k\geq 5$.

\end{proof}

\section{Conclusion}
The above discussion centers around two major themes: spectral
residues, and the distribution of coordinate values in a simplex of
bounded, ordered integer tuples.  In the first case, we have
demonstrated that the method of spectral residues leads to composition
sum identities with interesting interpretations.  We have considered
here parameterized recursive relations corresponding to first and
second-order linear differential equations in one independent
variable.  The next step in this line of inquiry would be to consider
other classes of parameterized recursive relations --- perhaps
non-linear, perhaps corresponding to partial differential equations
--- in the hope that new and useful composition sum identities would
follow.

In the second case, we have uncovered an interesting geometrical
property of the order simplex.  Theorems \ref{thrm:id1},
\ref{thrm:id2}, \ref{thrm:id3} support the conclusion that the middle
dimensions of an order simplex are more ``ample'' then the outer
dimensions.   However the 3 results we have been able to establish all
depend on very specific identities, and do not provide a general
tool for the investigation of this phenomenon.  To put it another way,
our results suggest the following
\begin{conjecture}
  Let $N$ be a natural number greater than $2$ and $d$ a natural
  number strictly less than $N/2-1$.  Let $n\geq N$ be another natural
  number.  For every sufficiently small cardinality $l$, there are
  more $l$-element subsets of $S^N(n)$ with distinct $x_{d+1}$
  coordinates, than there are $l$-element subsets with distinct $x_d$
  coordinates.
\end{conjecture}
\noindent It would also be interesting to see whether this conjecture holds if
we consider subsets of points with sparse, rather than distinct sets
of coordinate values.

Finally, Theorem \eqref{thrm:id2} deserves closer scrutiny, because it
describes a bijection of sets, rather than a mere comparison.  It is
tempting to conjecture that this bijection has an enumerative
explanation based on some combinatorial algorithm.


\begin{thebibliography}{XXXX}
  
\bibitem{Goulden} Goulden, I. and Jackson, D., {\em Combinatorial
    Enumeration}, Wiley, New York, 1983.

\bibitem{Hoggatt} V. E. Hoggat Jr., D. A. Lind, Compositions and
  Fibonacci numbers, {\em Fibonacci Quarterly}, {\bf 7} (1969),
  253--266.
  
\bibitem{Homenko} N. P. Homenko, V. V. Strok, Certain
  combinatorial identities for sums of composition coefficients. {\em
    Ukrain. Mat. Zh.} {\bf 23} (1971), 830--837.
  


\bibitem{Jordan} Jordan, C., {\em Calculus of Finite Differences},
  Chelsea, New York, 1947.

\bibitem{milson2}Milson, R., {\em Spectral residues of second-order
    differential equations: a new methodology for summation identities
    and inversion formulas}, preprint  math-ph/9912007.
  
\bibitem{Moser} Moser, L. and Whitney, E., Weighted Compositions, {\em
    Can. Math. Bull.}, 4:39--43 (1961).
  
\bibitem{wilf}Wilf, H., ``Generatingfunctionology'', Academic Press,
  1990.
\end{thebibliography}
\end{document}